\documentclass[12pt]{article}
\usepackage{graphpap}
\usepackage{stmaryrd}
\usepackage{amssymb}
\usepackage{amsmath}
\usepackage{algorithm}
\usepackage[noend]{algorithmic}

\def \qed {\hfill $\boxempty$}

\newtheorem{lem}{Lemma}
\newtheorem{Theorem}{Theorem}
\newtheorem{defi}{Definition}
\newtheorem{crl}{Corollary}
\newtheorem{prp}{Proposition}
\newtheorem{rmk}{Remark}
\newtheorem{xmp}{Example}
\newtheorem{clm}{Claim}
\newtheorem{op}{Problem}
\newtheorem{con}{Conjecture}

\def \bp {\begin{prp} \ }
\def \ep {\end{prp}}

\def \bc {\begin{crl} \ }
\def \ec {\end{crl}}

\def \bcon {\begin{con} \ }
\def \econ {\end{con}}

\def \thm {\begin{Theorem} \ }
\def \ethm {\end{Theorem}}

\def \bl {\begin{lem} \ }
\def \el {\end{lem}}

\def \bd {\begin{defi} \ \rm }
\def \ed {\end{defi}}

\def \brm {\begin{rmk} \ }
\def \erm {\end{rmk}}

\def \bxm {\begin{xmp} \ \rm }
\def \exm {\end{xmp}}

\def \bcm {\begin{clm} \ }
\def \ecm {\end{clm}}

\def \bop {\begin{op} \ }
\def \eop {\end{op}}

\def \nmr {\begin{enumerate}}
\def \enmr {\end{enumerate}}

\def \tmz {\begin{itemize}}
\def \etmz {\end{itemize}}

\def \nin {\noindent}
\def \bsk {\bigskip}

\def \msk {\medskip}
\def \pf {\nin{\bf Proof } \ }
\def \qed {\hfill $\Box$}

\def \es {\emptyset}

\begin{document}

\title {Domination game on forests\footnote{Research supported by the European Union and
Hungary and co-financed
 by the European Social Fund through the project T\'AMOP-4.2.2.C-11/1/KONV-2012-0004 - National Research Center
 for Development and Market Introduction of Advanced Information and Communication Technologies.}}

\author{Csilla Bujt\'{a}s\footnote{Email: bujtas@dcs.uni-pannon.hu }
\\
\normalsize
Department of Computer Science and Systems Technology \\
 \normalsize University of Pannonia\\
    \normalsize Veszpr\'{e}m,   Hungary
}
\date{}
\maketitle

\begin{abstract}
 In the domination game studied here, Dominator and Staller alternately choose a vertex of a graph
 $G$ and take it into a set $D$. The number of vertices dominated by the set $D$ must increase in each single
 turn and the game ends when $D$ becomes a dominating set of $G$. Dominator aims to minimize whilst Staller aims to maximize the
 number of turns (or equivalently, the size of the dominating set $D$ obtained at the
 end). Assuming that Dominator starts and both players play optimally, the
 number of turns is called the game domination number $\gamma_g(G)$ of $G$.

  Kinnersley, West and Zamani verified that $\gamma_g(G) \le 7n/11$ holds for every
  isolate-free $n$-vertex forest $G$ and they conjectured that the sharp upper bound is only $3n/5$.
  Here, we prove  the 3/5-conjecture for forests in which no two leaves are at distance 4 apart.
 Further, we establish an upper bound $\gamma_g(G) \le 5n/8$, which is valid for
 every isolate-free forest $G$.

\bigskip


\bigskip

\nin \textbf{AMS 2000 Subject Classification:}
 05C57,
 91A43,
 05C69

\end{abstract}

\vfil

\section{Introduction}

\subsection{Domination game}
 The  domination game  considered here
 was introduced in 2010 by Bre\v{s}ar,  Klav\v{z}ar and   Rall  \cite{BKR-SIAM},
  where the original idea is  attributed   to Henning (2003, personal
 communication).
  For this {\it domination game},  a graph $G$ is given  and two players, called  Dominator and
 Staller, take turns choosing a vertex  and taking it into a set $D$.
 Each vertex chosen dominates itself and its neighbors. The rule of
 the game prescribes that the set of vertices dominated by $D$ must
 be enlarged in each single turn. The game ends when no more legal
 moves can be made; that is, when $D$ becomes a dominating set of
 $G$.
 The goal of Dominator is to minimize, while that of Staller is to
 maximize the length of the game. Equivalently, Dominator wants a
 small dominating set $D$ and Staller wants $D$ to be as large as
 possible.
 The {\it game domination number} $\gamma_g(G)$ of $G$ is the
 number of turns in the game
 (equals the cardinality  of the dominating set $D$ obtained at the
 end) when Dominator starts the game and each of the two  players applies an optimal strategy.
 Analogously, the {\it Staller-start game domination number} $\gamma_g'(G)$ is the
 number of turns  when Staller begins   and the players   play
 optimally.

 \bigskip

 \subsection{Standard definitions}
 For a vertex $v\in V$ of a  graph $G=(V,E)$,  its {\em  open
  neighborhood}   is defined as   $N(v)=\{u : uv \in E\}$, whilst
  its {\em closed neighborhood} is   $N[v]=N(v) \cup \{v\}$.
  Then the {\em degree} $d(v)$ (or $d_G(v)$) of $v$ is just $|N(v)|$.
  Each vertex dominates itself and its neighbors, moreover
   a set $S\subseteq V$ dominates all vertices contained in
   $N[S]=\bigcup_{v\in S}N[v]$. A vertex set $D \subseteq V$ is
   called {\em dominating set} if $D$ dominates all vertices of $G$.
   The smallest cardinality of a dominating set $D$ is the {\em
   domination number}\/ $\gamma(G)$ of $G$.
   One can prove that
   $\gamma(G) \le \gamma_g(G) \le  2\gamma(G)-1$ and
   $\gamma(G) \le \gamma_g'(G) \le  2\gamma(G)$  hold.

In a tree, as usual,  a leaf is a vertex of degree 1, while a vertex
 having a leaf-neighbor is called {\em stem}.

 \bsk
 \subsection{Results on the domination game}

 The earlier papers discuss several aspects of the domination game,
 for example,  connections between $\gamma_g(G)$ and
 $\gamma_g'(G)$ \cite{BKR-SIAM, KWZ, Kos}, the game domination number of
 Cartesian products \cite{BKR-SIAM} moreover the difference between
 $\gamma_g(G)$ and $\gamma_g(H)$ when $H$ is a spanning subgraph of
 $G$ \cite{BKR}. The recent manuscript \cite{BDKK} discusses the
 possible changes of the game domination number when a vertex or an
 edge is deleted from the graph.

  From our point of view, the  following ``3/5-conjecture''
  and the related results
    are the most important ones.

  \begin{con}[Kinnersley,  West and Zamani, \cite{KWZ}]
  \label{conj1}
  If $G$ is an isolate-free forest of order $n$, then
  $$\gamma_g(G) \le \frac{3n}{5} \quad \mbox{and} \quad \gamma_g'(G) \le \frac{3n+2}{5}.$$
  \end{con}

  Conjecture \ref{conj1}    is proved to be true
  for graphs
   each of whose  components is
  a caterpillar\footnote{A caterpillar is a tree whose non-leaf
  vertices induce a path.} \cite{KWZ}.
  Additionally,  the authors of the   recent paper \cite{BKKR}
    identify all trees attaining this bound up to 20 vertices by
  computer search moreover construct  infinitely many trees
     satisfying $\gamma_g =3n/5$.
  As it follows, the bound $3n/5$ (if true) is sharp.

  \bigskip

  One of our contributions is the proof of Conjecture \ref{conj1} for the
  class of  forests in which no two leaves are connected by a path of length
  4. For this class of forests our upper bound $(3n+1)/5$ on
  $\gamma_g'$ is   slightly better than the bound conjectured in \cite{KWZ} for
   forests in general.

  \thm \label{th1}
  If $G$ is an isolate-free forest of order $n$ in which no two
  leaves have distance 4, then
  $$\gamma_g(G) \le \frac{3n}{5} \quad \mbox{and} \quad \gamma_g'(G) \le \frac{3n+1}{5}$$
  hold.
  \ethm

  Our proof, presented in Section \ref{sect3}, is based on a value-assignment  to the vertices,
  where the value of a vertex $v$ depends on
   the current status of $v$ in the game.
  Then, we describe a greedy-like strategy for Dominator  which
    ensures that the game ends within   $3n/5$ turns.
    We introduced this approach    in the conference paper
    \cite{JH}, where also Theorem~\ref{th1} was stated without a
    completely detailed proof.
     Now, the strategy described there is fine-tuned and the proof is extended by a
   a more detailed analysis to obtain a further result. This new general upper
   bound $5n/8$ concerns all isolate-free forests
   and
    improves the earlier bound
   $\gamma_g(G) \le  7n/11$, which was recently proved by
   Kinnersley,  West and Zamani  \cite{KWZ}.

 \thm \label{th2}
  If $G$ is an isolate-free forest of order $n$,
    then
  $$\gamma_g(G) \le \frac{5n}{8} \quad \mbox{and} \quad \gamma_g'(G) \le \frac{5n+2}{8}.$$
    \ethm

    \bsk

 The paper is organized as follows. In Section 2, the basic value-assignment
 is introduced and some general lemmas are obtained.
 Then, in Section 3, we describe the strategy and analyze the
 structure of the residual graph at some crucial points. In the last
 subsection of this part, we verify Theorems~\ref{th1} and \ref{th2}
 based on the previous lemmas. In Section 4, we make some concluding
 notes.



\section{Preliminaries} \label{sect2}

 At any moment of the game we have three different types of the
 vertices. We assign  them to different colors  and to different
  numbers of
 points. The letter $D$ always denotes the set of   vertices selected by the players up
 to the considered moment
  of the game. A vertex $v$ is dominated if   $v \in N[D]$,
   otherwise $v$
 is called undominated.

 \begin{itemize}
 \item A vertex   is {\it white}   and its value is 3 points if it is
 undominated.
 \item A vertex is {\it blue} and its value is 2 points if it is dominated
 but has at least one undominated neighbor.
 \item A vertex is {\it red} and its value is 0 point  if it and all of its
 neighbors are dominated.
 \end{itemize}

 Clearly,   selecting a red vertex would not enlarge the set of dominated
 vertices, hence this choice is not legal in the game. Also,
 selecting any vertex, the status of a red vertex will not change.
 Hence, red vertices can be ignored in the continuation of the
 game.
 On the other hand, blue vertices can be chosen later by any players as
 they have white neighbors, but edges connecting two blue vertices
 can be deleted.
 Therefore, at any moment of the game, graph $G$ will be  meant
 without  red vertices moreover without
  edges  joining two blue vertices.    This graph $G$ will be  called {\em residual graph}
  as it was introduced already in  \cite{KWZ}.
 Due to our definition, in a residual graph  each blue vertex has only  white neighbors and
 definitely has at least one.
 As relates white vertices, none of their neighbors and none of the
  edges incident with at least one white vertex were deleted.
  This implies the following statements.
  \bl \label{W-leaf}
  \tmz
  \item[$(i)$]
  If $v$ is a white vertex in a residual graph $G$, then $v$ has
  the same neighborhood in $G$ as it had at the beginning of the
  game. Particularly, if $v$ is a white leaf in $G$, then it was
  white leaf in each of the earlier residual graphs.
   \item[$(ii)$]
 If $G$ contains no isolated vertices at
 the beginning of the game, this property remains valid  for each
 residual graph throughout the game.
 \etmz
 \el

 When a vertex $v$ is played, it becomes red,   each white vertex
from $N(v)$ becomes either blue or red and additionally, each blue
leaf contained in $ N[N(v)]$ turns red. Further, if the game is
played on a tree, these are the only possible changes in colors.

\bsk

 The {\em value  $p(G)$ of a residual graph} $G$ is defined to be the sum of the values associated with its
 vertices. When a player selects
  a vertex, $p(G)$
 necessarily decreases. We say that the player gets (or seizes)  $q$ points in a turn if
 his move causes decrease $q$ in the value of $G$.

 Observe that in each turn
   the player either selects a white vertex which turns red
 (this means 3 points by itself, even without additional gain); or selects
 a blue vertex $v$ which turns red (2 points) moreover $v$ must have at least one white
 neighbor which becomes blue or red (at least 1 additional point).
 Hence, we have

 \begin{lem} \label{3-points}
 In each turn, the value of
 $G$   decreases by at least $3$ points.
 \end{lem}
 As a preparation  for proving our main theorems, we introduce some further
 notations and terminology.
 \tmz
 \item In general, at any moment of  the game, $G$   denotes the current residual
 graph.  However, if   preciseness requires, we also use the
 notation $G_i$ for the residual graph obtained after the $i$th turn
 of the game, moreover the graph given at the beginning is referred
 to
 as $G_0$.    Similarly, the number of white, blue and red vertices
 after the $i$th turn are denoted by $w_i$, $b_i$ and $r_i$,
 respectively, and we set $w_0=n$, $b_0=r_0=0$.
 Thus, $p(G_i)=3w_i+2b_i$ and the number of points  the player got
 in the $i$th turn is just the difference $p(G_{i-1})-p(G_i)$.
 Note that in the Dominator-start version, the $i$th turn belongs to Dominator, if $i$ is odd;
 otherwise it is Staller's turn.
 \item The subgraph of $G$ (or that of $G_i$) induced by the set of
 its white vertices is $G(W)$ (or  $(G_i(W)$, respectively).

 \item As relates colors, we use the abbreviations W, B and R.
 Hence, an R-vertex is a red vertex, a W-neighbor is a neighbor
 which is white and a B-leaf is a leaf of $G$ which is blue. Similarly,
 the notation $v:$ B$\rightarrow$R means that in the turn
 considered the color of $v$ changed from blue to red.
 Also,
 for a path subgraph of $G$, its type is denoted by the
 order of colors, for example BWB means a path on three vertices
 with the color-order indicated.

 \item A {\em critical $P_5$} is a path on five vertices whose both
 ends are  W-leaves and which is of type WWBWW. The unique blue vertex
 in a critical $P_5$ is called {\em critical center}.

 \etmz

At the end of this section we prove a further useful lemma.

  \bl \label{B-leaf}
 If the $i$th turn belongs to Dominator and the  residual graph $G_{i-1}$ contains a B-leaf in
 a component of order at least 3, then Dominator can seize at least
 7 points in the $i$th turn.
  \el
  \pf Assume that $G_{i-1}$ contains the B-leaf $v$ and
  let $u$ be its unique neighbor, which is definitely white.
  If $u$ has at least two W-neighbors, then Dominator   gets at
  least $7=2+3+1+1$ points by playing $u$.
  If $u$ has exactly one W-neighbor, say $u'$, then choosing $u'$,
  all the vertices $v$, $u$ and $u'$ become red, hence Dominator can
  seize at least $8=2+3+3$ points.
  If $u$ has no W-neighbor but the component consists of at least 3 vertices,
  then $u$ has a B-neighbor $v'$ which is
  different from $v$.
  In this case, if Dominator chooses  $v'$,  the value of $G_{i-1}$ decreases
  by at least
  $7=2+3+2$ points. \qed

 \bigskip

\section{Proof of the theorems} \label{sect3}

Here we prove our main results, Theorems \ref{th1} and \ref{th2}.
The two proofs are not separated, as they apply the same strategy
for Dominator, they
  proceed  by the same structural analysis and use  the same  lemmas.
 The special condition in Theorem \ref{th1}, namely the absence of leaves at distance four
 apart, will be used only in the final part of the proof.
 \bsk

 First, we consider the Dominator-start game on an isolate-free  $n$-vertex forest
 $G$,
 and describe a  strategy for Dominator
 which ensures the game to end within a limited number of turns whatever
 strategy is applied by Staller.
 In our presentation, the game is
 divided into four phases, some of which might be missing.
 For each Phase $i$ (for $i=1,2,3,4$) we give a strategy  prescribed for
 Dominator. Then,    Phase $i$ itself will be defined due to the
 applicability of
 the given strategies.
 \tmz
  \item  \underline{Strategy-Phase$(1)$} \quad In   his  turn,
  Dominator gets at least 7 points, moreover
    at least two vertices become red  in this turn.
  \item  \underline{Strategy-Phase$(2)$} \quad In   his  turn, Dominator
    gets at least 7 points.
    \item  \underline{Strategy-Phase$(3)$} \quad In   his  turn, Dominator
    gets at least 6 points.\\
    In Phase 3 we have two additional rules Dominator must apply:
    \tmz
    \item  \underline{$(R.3.1)$} \quad Dominator
    plays a vertex which results in the
 possible maximum gain achievable in that turn.
      \item  \underline{$(R.3.2)$} \quad Under the  rule $(R.3.1)$  Dominator prefers to
    play a W-stem having a W-leaf neighbor.
    \etmz
    \item  \underline{Strategy-Phase$(4)$} \quad In   his  turn, Dominator
    gets at least 3 points.
 \etmz
 Phase $i$ may start only  with the first turn of Dominator  when there is no
 applicable Strategy-Phase$(j)$ for any integers $1 \le j \le
 i-1 $. But it really  starts only if   Strategy-Phase$(i)$ can be applied in this
 turn, otherwise this phase is skipped.
 If Phase $i$ was not skipped, then it  ends just
 before the first turn of Dominator when Strategy-Phase$(i)$ is not applicable.
    Phase 1 is skipped only if
  all components of $G$ are of order 2.
  Let us emphasize that  we never go back to an earlier phase
 (no matter whether  it was ended or skipped). For example, at a point of  Phase 3, the changes in the
 structure of the residual graph might cause that Dominator can get 7 points,
   but then the game is continued in Phase 3.
   We remark  that,  by Lemma~\ref{3-points}, Dominator  always
 is able to get at least 3 points if the game is not   over yet.

  In general, we observe that each non-skipped phase begins with a
 turn of Dominator and ends with a turn of Staller with the only
 exception when Dominator ends the game and hence the current phase as
 well.

 \bsk
 To prove the theorems, we  will keep track of  the decrease in $p(G)$  from phase to phase, moreover
 analyze the structural properties of the residual graph at some
 points of the game.
 \bsk


 \subsection{Phase 1}
 Due to Strategy-Phase$(1)$ and Lemma \ref{3-points}, Dominator
 seizes
 at least 7 points and Staller gets at least 3 points in each of
 their Phase-1-turns. The extra  points seized above these limits
 are counted separately and will be put to use in Phase 3 when
 critical $P_5$-subgraphs are treated.
  Formally, for every $i \ge 1$ if the $i$th turn belongs to Phase
  1, we define
  $$ e_i= \left\{
 \begin{array}{cl}
  (p(G_{i-1})-p(G_i)-7 & \mbox{if $i$ is odd}\\
  (p(G_{i-1})-p(G_i)-3 & \mbox{if $i$ is even}
       \end{array}
 \right.
 $$
  Moreover, let $e^*=\sum_{i=1}^k e_i$, where $k$ is the number of turns belonging to
  Phase 1.
  As Dominator begins the phase, we have
  \bl \label{ph1}
  If Phase 1 consists of $k$ turns
   ($k\ge 0$), then the value of
 $G$   decreases by   $5k+e^*$ in this phase, where $e^* \ge
 0$.
  \el
 Further, we estimate  the number of critical
 centers.

 \bl \label{lemma5}
 Let Phase 1 consist  of $k$ turns and
  let $r_k$ denote  the number of
 red vertices at the end of Phase 1.
 Then, the number of vertices  which are critical centers in at least
 one  later residual graph $G_i$ ($i \ge k$)
  is at
 most $(r_k/3)+e^*$.
 \el
 \pf Consider the color-changes in the $i$th turn of Phase 1 and denote the number of
 vertices with changes    W$\rightarrow$R, B$\rightarrow$R, and
 W$\rightarrow$B
  by $x_1$, $x_2$ and $x_3$ respectively. Then, the change in the number of blue vertices
 is  $b_i-b_{i-1}=x_3-x_2$, and $p(G_{i-1})-p(G_i)=3x_1+2x_2+x_3$.
 \bsk

 First, assume that this is Dominator's turn. Then, Strategy-Phase$(1)$
 ensures that $r_{i}-r_{i-1} =x_1+x_2 \ge 2$ and hence
 $$e_i +1=3x_1+2x_2+x_3-6 = 3(x_1+x_2)-6 +x_3-x_2 \ge b_i-b_{i-1}.$$

 In the other case, when Staller moves, the vertex selected
 definitely
 becomes red and thus, $r_{i}-r_{i-1}=x_1+x_2 \ge 1$ holds. This implies
 $$e_i = 3x_1+2x_2+x_3-3 = 3(x_1+x_2)-3 +x_3-x_2 \ge b_i-b_{i-1}.$$

 Consequently, for any two consecutive moves in the phase
 $$e_i+e_{i+1} +1 \ge b_{i+1}-b_{i-1}  \qquad \mbox{and} \qquad
 r_{i+1}-r_{i-1} \ge 3$$
 hold.

 Note that if
   the game is  finished in Phase 1, the lemma clearly holds.
   Otherwise  $k$ is even, and
    for the number of blue and red vertices
  \begin{equation}
    \frac{k}{2}+ e^* \ge b_k \qquad \mbox{and} \qquad
   r_k \ge \frac{3k}{2}
   \end{equation}
   are valid, which yield
   \begin{equation}
   \frac{r_k}{3}+e^* \ge b_k.
   \end{equation}

   What remains to prove is that each vertex  which occurs as a critical center in any
   later residual graph  is already blue in $G_k$.
   Consider a  critical $P_5$ subgraph $v_1v_2v_3v_4v_5$   of
   a $G_i$ ($i \ge k$).
   As $v_1$ is a W-leaf in $G_i$, Lemma \ref{W-leaf} implies that it
   is also a W-leaf in $G_k$ at the end of Phase 1. Clearly, vertex $v_2$
   is white and $v_3$ is not red in $G_k$.
   Moreover, if $v_3$ was white in $G_k$,
   then Strategy-Phase$(1)$ could be applied in the $(k+1)$st turn, as
    Dominator could select $v_2$, which would cause the
        color-changes $v_1,v_2:$ W$\rightarrow$R and
    $v_3:$ W$\rightarrow$B. This cannot be the case as the $k$th
    turn finishes Phase 1. Therefore, each later critical center
    $v_3$ must be blue in $G_k$, and the lemma follows. \qed

\bsk

\subsection{Phase 2}

Our first statement is a direct consequence of the definition of
Phase 2, of Lemma~\ref{3-points} and of the fact that Dominator
starts the phase.

 \bl \label{Ph2}
  If Phase 2 consists of $k$ turns
   ($k\ge 0$), then the value of
 $G$ is decreased by at least  $5k$ in this phase.
  \el

  Our main observation concerning this phase is that the structure of
  the residual graph is quite restricted at the end of Phase 2.
 In a residual graph $G$, $v$ is a {\it single white vertex} (single-W) if it has only blue
   neighbors, that is, $N_{G(W)}(v)=\es$; and a {\it white pair} (W-pair) consists of  two
   W-vertices
   $u$ and $v$ for which $N_{G(W)}(u)=\{v\}$ and
     $N_{G(W)}(v)=\{u\}$ hold.

 \begin{lem} \label{Ph-2-end}
   At the end Phase 2 the residual graph $G$ has
   the following properties:
   \tmz
   \item[$(i)$] For each white vertex $v$,
   \tmz
   \item either $v$ is a single white vertex,
   \item or $v$  is in a   white pair.
   \etmz
   \item[$(ii)$] If a leaf   is contained in a component of order at
   least 3, then it is white.
   \item[$(iii)$] If a blue vertex $v$ belongs to a component of order at least 3
   and $v$ has a single white neighbor, then
      $v$ has exactly one further neighbor, which is necessarily
      from a white pair.
      \item[$(iv)$] Each blue vertex is of degree at most 4.
   \etmz
   Moreover, the above statements $(i)-(iv)$ are valid for every
   residual graph $G$ from which we have no possibility of choosing a
   vertex and attaining a gain of at least 7 points.
 \end{lem}

 \pf If the game is finished in Phase 2, then
 the residual graph in question contains no vertex, and there is
 nothing to prove. Otherwise, $G$ is a residual
 graph in which Strategy-Phase$(2)$ cannot be applied; that is, no
 choice of Dominator can cause a decrease of at least 7 points in the value of
 $G$.
   \tmz
  \item[$(i)$] Assume that $G(W)$ has a component of order at least
  3. Then, let  $v$ be one of the leaves of this
  component, and let $u$ be the only neighbor of $v$ in $G(W)$. As
  the component contains at least one further vertex, there exists a
  vertex $z$  for which  $z \in  N_{G(W)}(u)$ and $z \neq v$. Therefore,
  if Dominator plays vertex $u$, then $u$ and $v$ turn red and
  additionally $z$ becomes either blue or red. Consequently,
  Dominator could seize at least $7=3+3+1$ points. This cannot be
  the case, hence each component of $G(W)$ contains either one or
  two  vertices corresponding to the single-W vertices and to the
  W-pairs in $G$.
   \item[$(ii)$] Due to Lemma~\ref{B-leaf}, otherwise (when the leaf is blue) Dominator
   could get at least 7 points.
   \item[$(iii)$] Assume that $v$ is a B-vertex and has a single-W
   neighbor $u$. Then, choosing $v$, vertex $u$ becomes dominated
   and has no undominated neighbor. This already gives $2+3=5$
   points gain for Dominator. Due to $(ii)$, $v$ is not a leaf in
   $G$ and hence has a further W-neighbor $z$. If $z$ was a
   single-W vertex, the choice of $v$ would result in at least
   $5+3=8$ points, which contradicts our condition. Thus, $z$ is from a W-pair.
    Moreover, if $v$
   had three different W-neighbors, namely $u$, $z$ and $z'$, then
   Dominator could seize at least $7=5+1+1$ points by choosing $v$.
   As this is not the case, $v$ has exactly two neighbors $u$ and
   $z$ and the statement follows.
   \item[$(iv)$] If a blue vertex $v$ had five different neighbors,
   then all of them would be white and
     the choice of $v$ would give a gain of at least $7=2+5 \cdot 1$
   points.
  \etmz
  Finally, we observe that the same arguments are valid for any
  moment of the game, when Dominator has no possibility to get more
  than 6 points. \qed

  \bsk
  We remark that properties $(i)-(iii)$ were already satisfied by the residual graph at the end of
  Phase 1. This could be verified analogously to the above proof, but we do not do so,
  as we will not use this fact in the present paper.

\bsk

\subsection{Structural lemmas for later phases}

 Here, we prove some properties   which remain valid throughout
 Phases 3 and 4 (even if during Phase 3 Dominator has the possibility of
 seizing  7 or more points).

 \begin{lem} \label{Ph-2-cont}
   Throughout  Phases 3 and 4, each residual graph $G$ has
   the following properties:
   \tmz
   \item[$(i)$] Each white vertex  is either
     single white or    it is in a   white pair.
      \item[$(ii)$] If a blue vertex $v$ has a single white neighbor
      $u$, and $u$ has no blue-leaf neighbor, then $v$ has exactly
      one further neighbor $z$, which is either in a white pair or
      it is a single white vertex having a blue-leaf neighbor.
         \item[$(iii)$] Each  blue vertex is of degree at most 4.
   \etmz
 \end{lem}

 \pf
 \tmz
 \item[$(i)$]  As new white vertices do not arise during the game,
 this follows from Lemma~\ref{Ph-2-end}$(i)$.
 \item[$(ii)$]  If $u,u'$ formed a white pair at the end of Phase 2, then either it remains a white
 pair, or both $u$ and $u'$ turn red and are deleted, or one of them
 remains white and the other one becomes a B-leaf.
 By our conditions, $u$ is a single-W vertex without blue leaf in $G$,
 hence $u$ also was single-W at
 the end of Phase 2.
 Then, our
 statement follows from Lemma~\ref{Ph-2-end}$(iii)$.
 \item[$(iii)$] Consider a B-vertex $v$ of $G$. If $v$ was already blue at
 the end of Phase 2, then by Lemma~\ref{Ph-2-end}$(iv)$, it had at
 most four W-neighbors and hence, in any later residual graph $G$ its
 degree is at most 4. In the other case, when $v$ was a W-vertex at
 the end of Phase 2, it was either single-W, but then it could not
 be blue in $G$; or $v$ was in the W-pair $uv$ and now it is a
 B-leaf with the only neighbor $u$. \qed
 \etmz
 \bsk

 Applying   Lemma \ref{Ph-2-cont}, we  prove a further
 lemma, which says that if $G$ has  no component of type BWB, then
 the maximum achievable gain cannot be exactly 7 points.

  \begin{lem} \label{B-leaf2}
  In Phases 3 and 4, for any residual graph $G$ the
  following statements  hold.
  \tmz
  \item[$(i)$]
  If $G$ has a
  component which is not of the type BWB, but it is of order at least 3,
  moreover this component contains a  blue leaf,
  then Dominator can seize at least
  8 points.
  \item[$(ii)$] If there is no blue leaf in a component $C$ of order at least
  4, then Dominator cannot seize more than 6 points by playing a vertex
  from $C$.
  \etmz
  \end{lem}
  \pf
  \tmz
   \item[$(i)$]
  Consider a component satisfying the conditions of the lemma, a
  B-leaf $ v$ from it, and  the only W-neighbor $u$ of $v$.
 If $u$ is in the W-pair $uu'$,  Dominator can choose $u$ and then
 all the three vertices $v$, $u$ and $u'$ become red, and
 Dominator gets at least $8=2+3+3$ points.
  If $u$ is a single-W vertex and has at
  least three B-leaf neighbors (including vertex $v$), then the
  choice of $u$ results in a gain of at least $9=3+3\cdot 2$ points.
  If none of the previous cases holds, and as it is assumed, the component is not of the type
  BWB, then $u$ is a single-W vertex and
  has a B-neighbor $z$ which is not a leaf. If Dominator selects
  vertex $z$, then $v$, $u$ and $z$ become red, moreover at least one
  further white neighbor of $z$ turns   blue or red. Thus,
  Dominator seizes at least $8=2+3+2+1$ points.
   \item[$(ii)$] As there is no blue leaf in the component, selecting any vertex $v \in V(C)$ in a
   turn, only   $v$ and its W-neighbors  change their
   colors. We have the following cases due to the type of the vertex
   $v$ chosen.
   \tmz
   \item
   If $v$ is a  B-vertex which has no single white neighbor, then by
   Lemma~\ref{Ph-2-cont}$(iii)$, $v$ has at most four W-neighbors
   each of which turns blue. Therefore, Dominator may seize at most
   $6=2+4\cdot 1$ points, if he selects $v$.
   \item
   If $v$ is a B-vertex with a single-W neighbor $u$, then by
   Lemma~\ref{Ph-2-cont}$(ii)$ and by the absence of B-leaves, $v$ has
   exactly one further neighbor $z$ which is in a W-pair. Hence,
   selecting $v$  Dominator  gets exactly $6=2+3+1$ points.
   \item
   If $v$ is  single-W, then no vertex from $N(v)$ changes its
   color and therefore the gain is exactly 3 points.
   \item
   If $v$ is in the W-pair $vu$, the only color-changes are $v,u:$
   W$\rightarrow$R and hence, Dominator gets exactly 6 points.
  \qed
  \etmz \etmz
  \bsk

  \subsection{Phase 3: the crucial point}

  It was easy to see that the average decrease in the value $p(G)$ of the residual graph was
    at least 5 points per turn in the first two phases.
    We will see that this average   holds in the last phase.
    Also,
    if Staller gets at least 4 points in the $i$th turn of Phase 3, then
    together with the next turn of Dominator, when he seizes at least 6 points
      the desired average is
    attained locally. Hence, we focus on the turns when Staller gets only 3
    points.

    Recall that Strategy-Phase$(3)$ prescribes greedy selection for
    Dominator. Further, if he cannot get more than 6 points, the
    preferred choice is to dominate a (W-leaf,W-stem) pair.

  \begin{lem} \label{crucial}
  If Staller gets 3 points in the $i$th turn in Phase 3, then
  at least one of the following statements is true.
  \tmz
  \item[$(a)$]
   Dominator gets at least 8 points in the $(i-1)$st turn.
   \item[$(b)$]
  Dominator gets  at least 7 points in the $(i+1)$st
   turn.
   \item[$(c)$] Dominator chooses a white stem $v_2$  of a critical
   $P_5$ $v_1v_2v_3v_4v_5$ in the $(i-1)$st turn, and Staller
   selects the  center $v_3$ in the $i$th turn.
   \etmz
  \end{lem}
  \pf Assuming that Staller gets 3 points in the $i$th turn, we have two cases to consider.

   \paragraph{Case 1} Staller selects a single-W vertex $v$.\\
  As this choice results in only 3 points, $v$ is not from a
  component of order two, moreover $v$ has no B-leaf neighbor in
  $G_{i-1}$. Therefore, by Lemma~\ref{Ph-2-cont}$(ii)$, each
  $B$-neighbor $u$ of $v$ has exactly one further neighbor. Thus,
  after the selection of $v$ (that is, in $G_i$) $u$ is a B-leaf.
  Also, it follows from Lemma~\ref{Ph-2-cont}$(ii)$   that the
  component containing $u$ in $G_i$ is of order at least 3.
  Then, Lemma~\ref{B-leaf} and the greedy strategy of Dominator
   imply that $(b)$ holds.
    \paragraph{Case 2} Staller selects a B-vertex $v$.\\
    As he gets only 3 points, $N(v)=\{u\}$, where   $u$  is white but not a single-W
   vertex, moreover, $u$ has no B-leaf neighbor.
    Then, by Lemma~\ref{Ph-2-cont}$(i)$,
    $u$ must be  from a W-pair $uu'$, and
      after the move of Staller, $u$ becomes a B-leaf in $G_i$.

    If $u'$ is not a leaf in $G_i$, then the component of the B-leaf $u$ is
    of order at least 3, hence by Lemma~\ref{B-leaf} Dominator gets
    at least 7 points in the $(i+1)$st turn and $(b)$ holds.

    Suppose thus that $u'$ is a W-leaf in $G_i$ and hence,
    in $G_{i-1}$ and $G_{i-2}$, too. We also assume that $(a)$ is
    not valid, that is, Dominator could not get 8 or more points in
    the $(i-1)$st turn. Our goal is to prove that under these
    conditions $(c)$ is necessarily true.

    First, observe that $v$ was not a B-leaf in $G_{i-2}$
    (otherwise choosing $u$ Dominator could seize at least 8
    points).
    Similarly, $v$ was not a W-vertex in $G_{i-2}$, as this would
    mean a `W-triplet' in Phase 3.
    Consequently, $v$ was a non-leaf B-vertex in $G_{i-2}$ and had a
    further W-neighbor, say $z$.
    This component $C_{i-2}$ of $G_{i-2}$ contains $v$, $u$, $u'$ and $z$, hence
    its   order is at least 4.
    As   we assume   $(a)$ not to be valid,   by
    Lemma~\ref{B-leaf2}$(i)$ we can conclude that $C_{i-2}$ contains
    no B-leaf.
    Next, we apply Lemma~\ref{B-leaf2}$(ii)$ and obtain that
    Dominator can seize at most 6 points in the $(i-1)$st turn,
     and further,
    as Phase 3 is not finished at this time, he surely gets exactly
    6 points.
    Due to the rule $(R.3.2)$ given in Strategy-Phase$(3)$, he
    selects a W-stem of a W-leaf if there exists such a pair.
    Actually, there does exist one, as the pair $uu'$ is of this type.
    Since Dominator did not choose $u$ in the $(i-1)$st turn, he
    played another W-stem with a W-leaf. This caused   change in
    the color of $z$, so the only possibility is that Dominator
    selected the W-stem $z$ which had a W-leaf $z'$.

    Therefore, $z'zvuu'$ was a critical $P_5$ in $G_{i-2}$ and
    Dominator chose the W-stem $z$ in the $(i-1)$st turn and then
    Staller played the center in the $i$th turn. This satisfies
    $(c)$. \qed

    \bsk

    When case $(c)$ of Lemma~\ref{crucial} is realized in the game
    and neither $(a)$ nor $(b)$ holds,
    the $i$th turn  (when the center is selected) is called {\em critical turn}.
        Note that all of such turns belong to Phase 3.
     In the following lemma we estimate the number $c^*$ of critical
     turns.

     \bl \label{lemma11}
     Let $n_\ell$ denote the number of non-red vertices at the
     beginning of Phase 3, and let $c^*$ be the number of critical
     turns. Then, $5c^* \le n_\ell$ holds.
     \el
     \pf It is clear by  definition  that the $i$th turn might be
     critical only if Dominator's   choice in the $(i-1)$st turn
     and   Staller's choice in the next turn together change three vertices to be
     red in a component of order at least 5, moreover a new
     component of order 2 (of type BW) arises. These five vertices
     are associated with the $i$th critical turn.
     As they were non-red vertices at the beginning of the phase and no
     vertex is associated with more than one critical turn, the
     inequality follows. \qed

  \bsk

  \begin{lem} \label{ph3}
 If Phase 3 consists of $k$ turns ($k\ge 0$) and $c^*$ denotes the number of critical turns,
  then the value of
 $G$ has been decreased by at least $5k-c^*$ in this phase.
 \end{lem}
 \pf
   For the sake of simplicity, let the turns of Phase 3 be indexed
   from 1 to $k$, and $d_i$ (for $i$ odd) and $s_i$ (for $i$ even)
   denote the number of points Dominator or Staller seized in the
   $i$th turn, respectively. Hence, the value of the residual graph
   $G$ was decreased by
   $$P=\sum_{1 \le i \le k,       \mbox{\enskip \tiny{$i$  odd }}} d_i +
   \sum_{1 \le i \le k,      \mbox{\enskip \tiny{$i$  even }}} s_i.$$
   First, if $s_i=3$ and $d_{i-1}\ge 8$, we redefine $s_i=4$ and
   $d_{i-1}=7$.
   Then, if the $i$th turn of Phase 3 is critical, we
   increase $s_i$ from 3 to 4.
   For the sum $P'$ of the current values, the inequality $P' \le
   P+c^*$ holds.

   Now, consider the pairs $s_i+d_{i+1}$ where $i$ is even and $2
   \le i <k$.
   If $s_i=3$, neither $(a)$ nor $(c)$ from Lemma~\ref{crucial} is
   true for this turn, hence $(b)$ must be valid and
   $s_i+d_{i+1}\ge   10$ follows.
   If $s_i\ge 4$ and $i <k$, then   $d_{i+1}\ge 6$, and we have $s_i+d_{i+1}\ge
   10$ again.

   If $k$ is even and the last move of the phase is made by
   Staller, then $(b)$ from Lemma~\ref{crucial} cannot be true.
   Thus, $s_k \ge 4$ and $d_1+ s_k \ge 10$, from which
     $P' \ge 5k$.
     Similarly, if $k$ is odd, $P' \ge d_1 +10(k-1)/2 >5k$ holds,
     and the lemma follows.  \qed

     \bsk

     \subsection{Phase 4}

     We show that the structure of the residual graph is very simple
     throughout this phase.

     \bl \label{ph4}
       If Phase 3 consists of $k$ turns, then the value of
     $G$ has been decreased by exactly $5k$ in this phase.
     \el
      \pf Consider the residual graph $G$ which we have at the beginning
      of this phase. As Dominator cannot seize 6 or more points,
      there are no W-pairs, hence each W-vertex is single-W.
      Now, if a B-vertex $v$ had at least two neighbors, then
      selecting vertex $v$, all of its neighbors and also $v$ itself
      would turn red, and Dominator would seize at least 8 points.
      Therefore, each blue vertex is a leaf. It is also easy to see
      that each white vertex has no more than one B-leaf neighbor, and definitely has at least one,
      as there are no isolated vertices by Lemma~\ref{W-leaf}$(ii)$.

      Consequently, each component of $G$ is a
    $K_2$ with one white and one blue vertex. Therefore, no matter which
   vertex is selected, in each turn the value of the residual graph is decreased by
    exactly 5.    \qed

 \bsk

 \subsection{Finalizing the proofs}

 Here  we present the proofs of our theorems,
 based on the lemmas verified in the previous subsections.

\paragraph{Proof of Theorem 1} \enskip
 Consider an isolate-free forest $G$  in which no two leaves are at
 distance 4 apart.
 By Lemma~\ref{W-leaf}$(i)$,   no new white leaves arise. Thus, we have no
 critical $P_5$ subgraphs at any moment of the game, and there
 occur  no critical turns in Phase 3. 

 At the beginning, we have $p(G)=3n$ and this is
 decreased to zero during the game.
 By Lemmas~\ref{ph1}, \ref{Ph2}, \ref{ph3}, \ref{ph4} and by
 $c^*=0$,
  the average decrease in the value of the residual graph is at least
 5  points per turn
 for the Dominator-start game.
  Then, the desired  upper bound immediately follows:
 $$\gamma_g (G) \le \frac{p(G)}{5}= \frac{3n}{5}.$$

 For the Staller-start version, we may define a Phase 0 consisting of just the
 starting  turn indexed by 0. Recall that
  $G$ contains no isolated vertices and every vertex is white, which
  implies that in this   turn Staller gets at least 4 points. Then,
  our
  lemmas on the later phases remain valid and we have
  $$\gamma_g'(G) \le \frac{3n+1}{5}
   $$
   as stated. \qed

   \paragraph{Proof of Theorem 2} \enskip
   In this general case, we consider an isolate-free
   forest $G$ with $p(G)=3n$.
   If the described strategy yields a game with $t$
   turns, $e^*$ extra points in Phase 1 and $c^*$ critical turns in
   Phase 3, our
   Lemmas~\ref{ph1}, \ref{Ph2}, \ref{ph3} and \ref{ph4} imply
   $$t \le \frac{3n-e^*+c^*}{5}.$$
   The number $c^*$ of critical turns in Phase 3 cannot be greater
   than the number of critical centers at the beginning of this
   phase. Moreover, by Lemma~\ref{lemma5} the latter parameter is
   not greater than $(r_k/3)+e^*$, where $r_k$ is the number of the
   vertices turned red in Phase 1. Therefore,
   $$r_k \ge 3(c^*-e^*).$$
   On the other hand, by Lemma~\ref{lemma11}, for the number
   $n_\ell$ of vertices which are non-red at the beginning of Phase
   3,
   $$n_\ell \ge 5c^* \ge 5(c^*-e^*).$$
   Thus, we obtain
   $$n \ge r_k+n_\ell \ge 8(c^*-e^*)$$
   and then,
   $$\gamma_g(G) \le t \le \frac{3n+(n/8)}{5}=\frac{5n}{8}$$
   as stated in Theorem 2.

   Similarly to the proof of Theorem~\ref{th1}, the Staller-start game is
   treated by introducing Phase 0. As $G$ is isolate-free, Staller
   gets at least 4 points in the turn indexed by 0 and then, for the
   number of red and blue vertices $r_0\ge 1$ and $r_0+b_0 \ge 2$
   hold. Lemma~\ref{lemma5} for this Staller-start version must be
   modified, as in the   inequalities (1) and (2), parameters
   $r_k$  and $b_k$ must be replaced by $r_k-r_0$ and by $b_k-b_0$
   respectively. Otherwise, the proof proceeds in the same way.
   Thus, here we obtain

   \paragraph{Lemma 5'} \enskip {\it Let Phase 1 consist  of $k$ turns and
  let $r_k$ denote  the number of
 red vertices at the end of Phase 1.
 Then, the number of vertices  which are critical centers in at least
 one  later residual graph $G_i$ ($i \ge k$)
  is at
 most $((r_k-r_0)/3)+e^*+b_0$.}
 \msk

  Observe that the same upper bound holds for $c^*$.
 \msk

 Let us introduce the notation $e_0^*= 3r_0+b_0-5$, which is the number of extra points
 achieved
 above the desired average 5 points in the starting turn.
 Note that $e_0^*$ might equal $-1$, but otherwise it is non-negative.

 For the number $t'$ of turns in this game,
 $$t' \le
 \frac{3n-e^*-e^*_0+c^*}{5}=\frac{3n+1-[c^*-e^*-(e^*_0+1)]}{5}.$$
Applying Lemma 5',
$$3[c^*-e^*-(e^*_0+1)]\le r_k-10r_0+12 \le r_k+2,$$
and by Lemma~\ref{lemma11},
$$5[c^*-e^*-(e^*_0+1)]\le 5c^* \le n_\ell$$
is obtained. Then, we conclude the inequality
$$\gamma_g' \le t' \le
\frac{3n+1+\frac{n+2}{8}}{5}=\frac{5n+2}{8}$$ which proves the
theorem.   \qed

\bsk

\section{Concluding remarks}

Although Conjecture \ref{conj1} is a challenging open problem in
itself, we close this paper with the following  more general version
of the conjecture.

\begin{con}[Kinnersley,  West and Zamani, \cite{KWZ}]
   If $G$ is an isolate-free graph of order $n$, then
  $$\gamma_g(G) \le \frac{3n}{5} \quad \mbox{and} \quad \gamma_g'(G) \le \frac{3n+2}{5}.$$
  \econ

It is worth noting that a graph may have greater game domination
number than any of its spanning trees. Hence, even if an upper bound
on $\gamma_g$ is verified for forests, there is no trivial way to
conclude the same bound for graphs in general.

 The relation $\gamma_g(G)
\le \lceil 7n/10 \rceil$ is the best result, which has been
published up to now for this general case \cite{KWZ}.
 In the forthcoming manuscript \cite{CS2}, we will improve this upper bound significantly by using our proof
 technique, where we consider a greedy-type strategy under some
 value-assignment to the vertices.


\end{document}